\documentclass[11pt]{article}
\usepackage{amssymb, amsfonts, amsmath}
\usepackage{latexsym}

\providecommand{\MR}{\relax\ifhmode\unskip\space\fi MR }

\renewcommand{\leq}{\leqslant}
\renewcommand{\geq}{\geqslant}

\begin{document}

\bigskip

\title{The word problem for some uncountable groups\\ given by countable words}

\bigskip
\bigskip


\author{\hspace*{-20mm}O. Bogopolski\hspace*{50mm} A. Zastrow\\
\hspace*{0mm}\begin{minipage}[t]{150mm}
\small{Institute of Mathematics of Siberian}\hspace*{28mm} {\small Institute of mathematics}\\ \small{Branch of Russian Academy
of Sciences,}\hspace*{21.5mm} {\small University of Gdansk, Poland} \\ {\small Novosibirsk, Russia}\hspace*{54.2mm} \small{e-mail: zastrow@mat.ug.edu.pl}\\ {\small and D\"{u}sseldorf University, Germany} \\ \small{e-mail:
Oleg$\_$Bogopolski@yahoo.com}
\end{minipage}
\\}


\maketitle

\bigskip
\begin{abstract}

We investigate the fundamental group of Griffiths' space,
and the first singular homology group of this space and of the
Hawaiian Earring by using (countable) reduced tame words. We prove that
two such words represent the same element in the corresponding group
if and only if they can be carried to the same tame word by
a finite number of word transformations from a given list.
This enables us to construct elements with special
properties in these groups.
By applying this method we prove that the two homology groups
contain uncountably many different elements that can be represented by infinite concatenations
of countably many commutators of loops.
As another application we give a short proof that these homology groups contain the direct sum of $2^{\aleph_0}$ copies of $\mathbb{Q}$. Finally, we show that the fundamental group of Griffith's space contains $\mathbb{Q}$.

\end{abstract}


\date

\maketitle

\renewcommand{\leq}{\leqslant}
\renewcommand{\geq}{\geqslant}

\section{\bf Introduction}


According to the classical group theory, a group $G$ is generated by a subset $\varLambda$, if every element of $G$
is a {\it finite product} of elements of $\varLambda\cup \varLambda^{-1}$.
However, for some topological groups, it is more natural to describe their elements by using (appropriately defined) infinite products of their ``basic'' elements.

Consider a one-point union of countably many circles $p_n$, $n\in \mathbb{N}$.
If we introduce the path metric on this set by saying that the length of every $p_n$ is 1,
we get a topological space whose fundamental group is the free group with the
basis $P=\{p_n\,|\, n\in \mathbb{N}\}$.
If we set the length of every $p_n$ to be $\frac{1}{n}$,
we get another topological space, the Hawaiian Earring $Z$.
Its fundamental group
is not free (see a short proof in~\cite{S}), but it can be well described by using special countable words over the alphabet $P$
(see~\cite{CC1}, for instance).

Another interesting example gives Griffith's space $Y$.
This remarkable space was constructed by H.B.~Griffiths in the fifties (see \cite[page 185]{Gr1}) as a one-point union of two cones over Hawaiian Earrings (see Figure~3 and Definition~3.3),
for showing that a one-point union of two spaces
with trivial fundamental groups can have a non-trivial and even uncountable fundamental
group. The fundamental group of Griffith's space can be also described by using special countable words over the alphabet $P\cup Q$, where $P$ is as above and $Q=\{q_n\,|\, n\in \mathbb{N}\}$ is a copy of $P$ (see Theorem~3.4).

In this paper we investigate the word problem in the first singular homology group and in the fundamental group
of the Hawaiian Earring and of Griffith's space.
Since all these groups are uncountable, an algorithmic solution is impossible.
We suggest another kind of solution by finding a finite number of transformations of countable words,
such that if two countable words $w_1,w_2$ represent the same element of a group, then one can get $w_2$ from $w_1$ by using a finite number of such transformations.

To be more precise, we explain
the second of the main Theorems~8.1--8.3.
Let $P=\{p_1,p_2,\dots
\}$ and $Q=\{q_1,q_2,\dots \}$ be two disjoint infinite countable alphabets.
In Section 2 we define the group $\mathcal{W}(P)$, consisting of all reduced countable words over the alphabet $P$, where every letter of $P\cup P^{-1}$ occurs only finitely many times.
This group can be identified with the group $\mathcal{P}$ introduced by Higman in~\cite[page 80]{H}, and  by
the result of Morgan and Morrison in~\cite[Theorem 4.1]{MM} it is canonically isomorphic to
$\pi_1(Z)$.
Thus, for the Hawaiian Earring~$Z$ we have the canonical homomorphisms
$$\mathcal{W}(P)\overset{\cong}{\longrightarrow}\pi_1(Z) \longrightarrow H_1(Z).$$
For Griffiths' space $Y$ we have the canonical homomorphisms
$$\mathcal{W}(P\cup Q)\longrightarrow \frac{\mathcal{W}(P\cup Q)}{\langle\!\langle \mathcal{W}(P),\mathcal{W}(Q)\rangle\!\rangle}
\overset{\cong}{\longrightarrow}\pi_1(Y) \longrightarrow H_1(Y),$$
where $\langle\!\langle \rangle\!\rangle$ means the normal closure, see Theorem~3.4 for details.

Theorem~8.2 states that two countable words from $\mathcal{W}(P\cup Q)$ represent the
same element of $H_1(Y)$ if and only if  each of them can be transformed to the same countable word
by applying a finite number of transformations of type {\rm (i)}, then of  type {\rm (ii)} and finally of type~{\rm (iii)}:

\medskip

(i)\hspace*{2mm}  Deletion of a countable subword, containing only $p$- or $q$-letters;

(ii)\hspace*{1mm}  Deletion of two distinct mutually inverse countable subwords: $AXBX^{-1}C\rightarrow ABC$;

(iii) Permuting two consecutive countable subwords.

\medskip

These transformations come from transformations of some associated 2-complexes, which we call arch-line-band systems; we define them in Section 4. These systems are closely related, but do not coincide
with the band systems introduced in the paper of M.~Bestvina and M.~Feighn \cite{BF}
for studying the stable actions of groups on real trees; see also the arc-surface technique
in the proof of J.W. Cannon and G.R. Conner~\cite[Theorem~4.3]{CC1}.

Theorems~8.1--8.3 enable us to construct and compare elements of the above groups with prescribed properties.
We deduce two corollaries.

Corollary~8.7 states that $H_1(Z)$ and $H_1(Y)$ contain uncountably many elements,
which can be represented as infinite products of  commutators.
This contrasts with the fact that any finite product of commutators
in the abelianization of a group is trivial.
In the forthcoming paper~\cite{BZ}, we investigate Karimov's space $K$ (see~\cite{Ka}) and make this contrast as
strong as possible. We show, that $H_1(K)$ is uncountable, and that each element
of $H_1(K)$ can be represented as an infinite commutator product.

Note that in~\cite[page~76]{H} Higman proved, in fact, that $H_1(Z)$ contains a nonzero element,
which can be represented as an infinite commutator product.
His short and nice proof uses a representation of a free group
as a group of formal power series in non-commuting variables, that was invented by W.~Magnus in~\cite{M}.

Corollary~9.2 contains a rather unexpected claim that the fundamental group of Griffith's space contains a subgroup isomorphic to $\mathbb{Q}$. The same fact for the groups $H_1(Z)$ and $H_1(Y)$ was known earlier, see~\cite{Eda3,Eda2}.
In the proof we use extensively a combinatoric of special countable words.

Writing this paper, we have realized a universal role of groups $\mathcal{W}(\varLambda)$ for (uncountable)
combinatorial group theory.
Using these groups instead of free groups, O.~Bogopolski and W.~Singhof developed in~\cite{BogSin} a theory of generalized presentations of groups and gave generalized presentations of arbitrary permutation group $\varSigma(X)$ and of the automorphism group of the free group of infinite countable rang, ${\text{\rm Aut}}(F_{\omega})$. A different approach to generalized presentations of groups was suggested by A.J.~Sieradski in~\cite{Sier}; sie also a paper and a preprint of W.A.~Bogley and A.J.~Sieradski~\cite{BS1,BS4}.


\bigskip

\section{\bf Countable words}

\bigskip

A theory of infinite words was independently developed by J.W.~Cannon and G.R.~Conner~\cite{CC1}, K.~Eda~\cite{Eda1,Eda4}, and A.~Zastrow~\cite{Z1,Z2,Z3}. Since all these approaches slightly differ from each other, we give here
an account of this theory in the form, which is appropriate for our aims.

\medskip

{\bf Definition 2.1.} Let $X$ be an alphabet,
$X^{-1}=\{x^{-1}\,|\, x\in X\}$ and
$X^{\pm}=X\cup X^{-1}$.
A {\it countable word} over $X$ is  sequence $(x_i)_{i\in I}$, where $I$ is an arbitrary linearly ordered
countable set and $x_i\in X^{\pm}$. The order on $I$ will be denoted by $\preccurlyeq$.
E.g., the natural numbers $\mathbb{ N}$ with the classical order
and $\mathbb N\times \mathbb N$ with lexicographical order are different types of countable linearly ordered
sets.
There are uncountably many
types of countable linear orders, e.g.\ all ordinal numbers that are smaller than the smallest
uncountable ordinal. Two countable words $U=(x_i)_{i\in I}$ and $V=(x_j)_{j\in J}$ are called {\it equal} if there
exists an order preserving bijection $\varphi:I\rightarrow J$, such that $x_{\varphi(i)}=x_i$ for all $i\in I$.
In this case we write $U\equiv V$ .

A subset $J$ of $I$ is called {\it connected}
if,   for any two elements $j_1,j_2\in J$ and any element $i\in I$, we have
$$(j_1\preccurlyeq i\preccurlyeq j_2)\Rightarrow  (i\in J).$$

Let $W=(x_i)_{i\in I}$ be a countable word over $X$.
A {\it subword} of $W$ is a sequence $(x_j)_{j\in J}$, where $J$ is a connected subset of $I$ with the induced order.
Two subwords $U=(x_p)_{p\in P}$ and $V=(x_q)_{q\in Q}$ of $w$ are called {\it consecutive} if $P\cup Q$ is a closed subset of $I$ and $P\cap Q=\emptyset$.

If we have two countable words $W_1=(x_i)_{i\in I_1}$ and $W_2=(x_i)_{i\in I_2}$, we can {\it concatenate} them,
i.e.\  build the new countable word $W=(x_i)_{i\in I_1\cup I_2}$, where we assume that the elements of $I_1$ are smaller than the elements of $I_2$. In this case we write $W=W_1W_2$, not using the dot between $W_1$ and $W_2$.

In order to define (ir)reducible countable words we first need:

\medskip


{\bf Definition 2.2.} (a visualization of countable words)

To each countable word $W=(x_i)_{i\in I}$ over the alphabet $X$ we associate

\hspace*{5mm} (i) a closed segment $[a,b]$ on the real axis of the Euclidean plane together with

\hspace*{5mm} (ii) a set of points in the open interval $(a,b)$, each
marked by a letter from $X^{\pm}$ so, that moving from $a$ to $b$ along this segment we read the word $W$.

The segment together with the set of its marked points will be called {\it $W$-segment}.
For brevity we will call these points {\it letters}.
Note, that the visualization is not unique.  The following notion will serve to visualize a cancelation process
in a countable word.

\medskip

{\bf Definition 2.3.} An {\it arch-system based on the $W$-segment} is a subset of the plane, consisting of the $W$-segment and of arches (i.e.\ half-circles) which are  attached to this segment
according to the following rules:

\hspace*{5mm} (1) all arches lie in the upper half-plane with respect to the real axis;

\hspace*{5mm} (2) each arch connects a letter with one of its inverses;

\hspace*{5mm} (3) different arches do not intersect;


\hspace*{5mm} (4) if two letters are connected by an arch, then
each letter between them is connected by an arch with another letter between them.

\medskip

\input{Cancel1.pic}

\begin{center}
Fig.~1
\end{center}

An arch-system based on the $W$-segment is called {\it complete}, if each letter of $W$ belongs to an arch of this system.
The endpoints of an arch $a$ are denoted by $\alpha(a)$ and $\omega(a)$; we
assume that $\alpha(a)<\omega(a)$ and call them the {\it initial} and the {\it terminal} points of $a$, respectively.



Two arches $a,b$ of a given arch-system are called {\it parallel} if after possible replacement of $a$ by $b$ and $b$ by $a$ the following holds:

1) $\alpha(a)<\alpha(b)<\omega(b)< \omega(a)$;

2) for any arch $c$ of the arch-system, if $\alpha(a)<\alpha(c)<\alpha(b)$ then $\omega(b)<\omega(c)<\omega(a)$.


\medskip

{\bf Definition 2.4.} A countable word $W$ is called {\it irreducible} (or {\it reduced}) if there does not exist an arch-system based on the $W$-segment with a nonempty set of arches. A countable word $W$ is called {\it reducible} if it is not irreducible.

\medskip

It is immediate from the definition that subwords of reduced words are reduced.




Let $U=(x_i)_{i\in I}$ be a countable word. A word $V=(x_j)_{j\in J}$ is called {\it inverse} to $U$ and is denoted by $U^{-1}$, if there exists
an order-reversing bijection $\varphi: I\rightarrow J$ such that
$x_{\varphi(i)}=x_i^{-1}$ for any $i\in I$. 
With $I=\varnothing$ the {\it``empty word''}\/ is also defined.


\medskip

Note, that reducible countable words do not necessarily contain two
consecutive inverse subwords, see an example in \cite[page 227]{CC1}.

In the following we discuss, how to define the ``reduced form'' of a countable word $W$, or in other words: how ``to reduce'' a countable word.
One of the natural definitions could be the following: we delete from $W$ the letters corresponding to the endpoints of arches of a maximal arch-system based on the $W$-segment and we will call the remaining countable word the reduced form of $W$. However different maximal arch-systems can give different reduced forms. The easiest example of this phenomenon is the following:


$$\begin{array}{ll}x\,\overset{\frown}{x^{-1}x}\,\overset{\frown}{x^{-1}x}\, \overset{\!\!\!\!\!\!\!\frown}{x^{-1}\dots} \rightarrow x;\vspace*{2mm} \\ \overset{\!\!\!\!\!\frown}{x\,x^{-1}}\overset{\!\!\!\!\!\frown}{x\,x^{-1}}\overset{\!\!\!\!\!\frown}{x\,x^{-1}}\dots \rightarrow \varnothing.
\end{array}$$
Moreover, the concatenation of two reduced countable words may have different reduced forms. Indeed, if $U=xxx\dots$ is infinite,
then $U$ and $U^{-1}$ are two reduced countable words, but $UU^{-1}$ can be reduced to $x^n$ for any integer $n$.
To avoid such unlikable effects, we will work in a smaller class of countable words
(that will be sufficient for our goals).




\medskip

 {\bf Definition 2.5.} A countable word $W=(x_i)_{i\in I}$ is called {\it restricted} if every letter $x\in X^{\pm}$ occurs only finitely many times in $W$.



\medskip


{\bf Definition 2.6.} Let $W$ be a restricted countable word.
A word $U$ is called the {\it reduced form of} $W$ if $U$ can be obtained from $W$
by deleting the letters corresponding to the endpoints of a maximal arch-system for $W$.

\medskip
Clearly, for finite words we get the usual reduced form, which is, of course, unique.
In Theorem~2.9, we will show that the reduced form of a restricted countable word is unique.
This theorem is not new (see \cite[Theorem~3.9]{CC1} of J.W.~Cannon and G.R.~Conner), however for completeness we give a proof. In \cite{CC1}, it was proved,
using topological properties of the plane, namely with the help of Jordan's theorem.
Our proof is purely algebraic and it was inspired by the paper \cite{Eda1} of Eda.




\medskip

For any subset $F\subseteq X$, let $W_F$ denote
the word obtained from $W$ by deleting all letters of $(X\setminus F)^{\pm}$. Similarly, if $\mathcal{A}$ is an arch-system for $W$, we denote by $\mathcal{A}_F$ the system that is obtained from $\mathcal{A}$ by deleting all arches with endpoints in $(X\setminus F)^{\pm}$. Clearly, $\mathcal{A}_F$ is an arch-system for $W_F$.
For any finite word $W$ we denote by $[W]$ its reduced form. 

\medskip

{\bf Lemma 2.7.} {\it Let $U$, $V$ be two  reduced forms
of a restricted countable word $W$ over $X$. Then for any finite
$F\subseteq X$ holds $[U_F]\equiv [V_F]$.} 

\medskip

{\it Proof.} Let $U$ be obtained from $W$ by deleting the
letter-endpoints of a maximal arch-system~$\mathcal{A}$. Then, $[U_F]$ is obtained from $W$ by the following consecutive operations:

\medskip

1) Remove letter-endpoints of all arches of $\mathcal{A}$.

2) Remove all letters of $(X\setminus F)^{\pm}$.

3) Reduce the resulting word in the free group with base $F$.

\medskip

This is equivalent to

\medskip

1) Remove all letters of $(X\setminus F)^{\pm}$.

2) Remove letter-endpoints of all arches of $\mathcal{A}_F$.

3) Reduce the resulting word in the free group with base $F$.

\medskip

This is equivalent to

\medskip

1) Remove all letters of $(X\setminus F)^{\pm}$.

2) Reduce the resulting word in the free group with base $F$.

\medskip

Similarly $[V_F]$ can be obtained from $W$ by the last two operations. This completes the proof. \hfill $\Box$

\medskip

Note that a finite word is reducible to the empty word if and only if there is a complete arch-system for it.

\medskip

{\bf Lemma 2.8.} {\it Let $W$ be a restricted countable word over $X$. If for any finite subset $F\subseteq X$,
there is a complete arch-system for $W_F$,
then there is a complete arch-system for $W$.}

\medskip

{\it Proof.} We may assume that $X$ is countable, say $X=\{ x_1,x_2,\dots \}$.
Denote $X_k=\{x_1,\dots ,x_k\}$, $k=1,2\dots$. For each $k\in \mathbb{N}$ we choose a complete arch-system associated with $W_{X_k}$ and denote it by $\mathcal{A}(W_{X_k})$. The restrictions of $\mathcal{A}(W_{X_k})$ to the subword $W_{\{x_1\}}$, where $k=1,2,\dots $, are complete arch-systems for $W_{\{x_1\}}$. Since there is only a finite number of arch-systems for the finite word $W_{\{x_1\}}$, we can choose an infinite subsequence $\mathcal{A}(W_{X_{k_1}})$, $\mathcal{A}(W_{X_{k_2}}),\dots $ with the same restrictions to $W_{\{x_1\}}$; denote the common restriction by $\mathcal{A}_1$. From this sequence, we can choose an infinite subsequence with the same restrictions to $W_{\{x_2\}}$; denote the common restriction  by $\mathcal{A}_2$, and so on. The systems $\mathcal{A}_1, \mathcal{A}_2,
\dots$ do not conflict (i.e. their arcs do not intersect each other), since any two of them are contained in a common finite arch-system. Therefore their union forms a complete arch-system for $W$. \hfill $\Box$

\bigskip

The idea of the proof of the following proposition is taken from~\cite[Theorem~1.4]{Eda4}.

\medskip

{\bf Proposition 2.9.} {\it Each restricted countable word $W$ over $X$ has a unique reduced form.}

\medskip

{\it Proof.} The existence follows straightforward by applying Zorn's lemma to the
set of all arch-systems based on $W$, which is partially ordered by inclusion.

Let us prove the uniqueness. Let $U=(x_i)_{i\in I}$, $V=(x_j)_{j\in J}$ be two reduced forms of $W$. By Lemma 2.7, for any finite $F\subseteq X$ holds $[U_F]\equiv [V_F]$.
Let $\alpha\in X$. Consider $U$ written in the form $U\equiv U_0A_1U_1A_2U_2\dots A_k U_k$,
where each $A_i$ is a positive or negative power of $\alpha$ and $U_i$ does not contain $\alpha$ and $\alpha^{-1}$, and $U_i$ is nonempty for $i\neq 0,k$.
Since $U$ is reduced, all $U_i$'s are also reduced and so, by Lemma~2.8, there exist finite subsets $F_i\subseteq X$ ($i=1,\dots ,k-1$), such that $[(U_i)_{F_i}]\neq \varnothing$. Hence for any finite subset $F\subseteq X$ containing $\cup F_i\cup \{\alpha\}$, we have $[(U_i)_{F}]\neq \varnothing$ for all $i\neq 0,k$, and so $[U_F]\equiv [(U_0)_{F}]A_1 [(U_1)_{F}]A_2 \dots A_k[(U_k)_{F}]$. In particular, $U_{\{\alpha\}}\equiv [U_F]_{\{\alpha\}}$.
Similarly, $V_{\{\alpha\}}\equiv[V_G]_{\{\alpha\}}$ for some finite subset $G\subseteq X$.

By Lemma 2.7, $[U_E]\equiv [V_E]$ for all finite subsets $E\subseteq X$. Hence, for a sufficiently large finite subset $E\subseteq X$ we have
$$U_{\{\alpha\}}\equiv [U_E]_{\{\alpha\}}\equiv [V_E]_{\{\alpha\}}\equiv V_{\{\alpha\}}.\eqno{(\dag)}$$

In the same way (defining the subwords $A_i$ now as all maximal subwords that contain only the letters $\alpha^{\pm 1}$ and $\beta^{\pm 1}$) one can prove that for any two letters $\alpha,\beta\in X$, there exists a finite subset $L\subseteq X$, such that
$$U_{\{\alpha,\beta\}}\equiv [U_L]_{\{\alpha,\beta\}}\equiv [V_L]_{\{\alpha,\beta\}}\equiv V_{\{\alpha,\beta\}}.\eqno{(\ddag)}$$

\bigskip

Formula $(\dag)$ enables to construct a bijection $\varphi:I\rightarrow J$, such that $x_{\varphi(i)}=x_i$ for all $i\in I$.
Formula $(\ddag)$ shows that $\varphi$ respects the orders on $I$ and $J$. So, we have $U\equiv V$. \hfill $\Box$

\medskip

%


{\bf Definition 2.10.} Let $W_1,W_2$ be two restricted reduced countable words over $X$. Their {\it product}, denoted $W_1\cdot W_2$, is the reduced form of the concatenation $W_1W_2$. By Definition~2.4 and Theorem~2.9, the set of all reduced countable words over $X$ with respect to this multiplication is a group. We denote this group by $\mathcal{W}(X)$.

\medskip

We stress, that $\mathcal{W}(X)$ is not generated by $X$, if $X$ is an infinite set. In Section~3 this group will be related to the fundamental groups of the Hawaiian Earrings and of Griffiths' space.
For technical reasons we will work with tame words, which we define next.


\medskip

{\bf Definition 2.11.} A restricted reduced countable word is called a {\it reduced}
(or {\it irreducible}) {\it tame word}. A finite concatenation of such words is called a {\it tame word}.

\medskip

{\bf Remark 2.12} {The class of tame words over $X$ is closed under taking subwords and finite concatenations}.

\medskip

In general, a tame word $W$ can be presented as a concatenation of reduced tame words in different ways.
The following lemma states, that at least one of these presentations
can be transformed to the reduced form of $W$
in the same way as in a free group.

\medskip

{\bf Lemma 2.13.} {\it Any tame word $W$ can be presented as a concatenation of nonenempty reduced tame words $W\equiv v_1v_2\dots v_k$, $k\geqslant 0$, such that the following process yields the reduced form of~$W$:

\vspace*{-7mm}
\noindent $(\dag)$
\hspace*{0.5cm}{\begin{minipage}[]{15cm}{\vspace*{10mm}
if some adjacent factors $v_i,v_{i+1}$ are mutually inverse, we delete them,
and if this leaves us with mutually inverse adjacent factors, we also delete them.
We repeat these steps until none of the remaining adjacent factors are mutually inverse.}
\end{minipage}}}






\bigskip

{\it Proof.} As $W$ is a tame word, $W$ is a finite concatenation of reduced tame words:
$W\equiv u_1u_2\dots u_l$. We may assume that all $u_i$ are nonempty and that $l\geqslant 2$, and we will use induction by~$l$.

Let $l=2$, that is $W\equiv u_1u_2$.
If $W$ is reduced, we can set $v_1=u_1$ and $v_2=u_2$.
If $W$ is not reduced, there exists a nonempty maximal arch-system $\mathcal{A}$ for $W$; since $u_1,u_2$ are reduced,
no arch of $\mathcal{A}$ has both endpoints in $u_1$ or in $u_2$.
Hence the initial points of these arches lie in $u_1$ and the terminal one in $u_2$, and all the arches are parallel.
Therefore there exist concatenations $u_1\equiv v_1v_2$, $u_2\equiv v_3v_4$, such that $v_2,v_3$ are mutually inverse and $v_1v_4$ is reduced. Thus, we can set $W\equiv v_1v_2v_3v_4\equiv v_1v_2v_2^{-1}v_4$ and the above described process yields $v_1v_4$ as the reduced form of~$W$.

Now let $l\geqslant 3$. By induction, $W'\equiv u_1u_2\dots u_{l-1}$ can be presented as a concatenation $W'\equiv w_1w_2\dots w_k$, with the reduced form $W'_{\text{\rm red}}\equiv w_{i_1}w_{i_2}\dots w_{i_s}$, obtained from the concatenation by applying the process $(\dag)$. Also  $W'_{\text{\rm red}}u_{l}$ (as in case $l=2$),
can be presented as $v_1v_2v_3v_4$, where $W'_{\text{\rm red}}\equiv v_1v_2$, $u_l\equiv v_3v_4$, and $v_2,v_3$ are mutually inverse, and $v_1v_4$ is reduced.

Then there exists $1\leqslant j\leqslant s$ and a subdivision $w_{i_j}\equiv pq$ such that $v_2\equiv qw_{i_{j+1}}\dots w_{i_{s-1}}w_{i_s}$,
and so $v_3\equiv w_{i_s}^{-1}w_{i_{s-1}}^{-1}\dots w_{i_{j+1}}^{-1}q^{-1}$.
The desired concatenation is $$W\equiv w_1w_2\dots w_{i_{j}-1}pqw_{i_{j}+1}\dots w_{k-1}w_kw_{i_s}^{-1}w_{i_{s-1}}^{-1}\dots w_{i_{j+1}}^{-1}q^{-1}v_4.$$ \hfill $\Box$




\section{\bf Hawaiian Earring and Griffiths' space}

{\bf Definition 3.1.} The {\it Hawaiian Earring} is the topological space, which is the (countable) union of circles of radii $\frac{1}{n}$, $n=1,2,\dots $, embedded into the Euclidean plane in such a way that they have only
one common point (called the {\it base-point}).
The topology of the Hawaiian Earring is induced by the topology of the plane.

\medskip

\input{HawaiiEarr3.pic}

\begin{center}
Fig.~2
\end{center}

Usually one assumes that these circles lie on same side of the common tangent (see Fig.~2).
If we name the circle of radius ${1\over n}$ by $p_n$, we will denote our Hawaiian Earring by $H(p_1,p_2,\dots)$ or shortly $H(P)$, where $P=\{p_1,p_2\dots \}$. For convenience we will use the same name $p_n$
for the path traversing $p_n$ in counterclockwise direction, and also for the homotopy class of this path.

Note, that each homotopy class in $\pi_1(H(P))$ can be represented by a path without backtracking inside of each $p_i$. Because of continuity the path cannot pass through $p_i$ infinitely many times. This yields
the natural map $\pi_1(H(P))\rightarrow \mathcal{W}(P)$. The following theorem is well known; see~\cite[Theorem~4.1]{MM}.

\medskip

{\bf Theorem 3.2.} {\it The fundamental group of the Hawaiian Earring is isomorphic
to the group $\mathcal{W}(P)$.}

\bigskip







The Hawaiian Earring has the following remarkable properties, some of them are not shared by tame spaces (e.g. manifolds or triangulable spaces).

\medskip

$\bullet$
It is not semi-locally simply connected.\\
(A space $X$ is called {\it semi-locally simply connected} if every point $x\in X$ has a neighborhood $N$ such that
the homomorphims $\pi_1(N,x)\mapsto \pi_1(X,x)$, induced by the inclusion map of $N$ in $X$, is trivial.)

$\bullet$ It does not admit a universal covering space.\\
(For connected and locally path connected spaces this is
equivalent to semi-local simple connectedness, see Corollary~14 in Ch.~2, Sec.~5 of~\cite{Sp}; see also~\cite[Example 2.5.17]{Sp}.)


$\bullet$
It is compact, but its fundamental group is uncountable.\\
(This is an easy exercise, see for example \cite[Theorem 2.5 (1)]{CC1}.)

$\bullet$
It is one-dimensional, but its fundamental group is not free.\\
(This follows from \cite[Theorem~4.1]{MM} and a remark in \cite[page 80]{H}. See also a short proof in \cite{S}.)

$\bullet$
Its fundamental group contains countable subgroups, which are not free.\\
(See \cite[page 80]{H} or \cite[Theorem 3.9 (iii)]{Z3}.)

$\bullet$ Every finitely generated subgroup of its fundamental group is free.\\
(Since the fundamental group embeds in an inverse limit of free groups, see \cite[Theorem~4.1]{MM}.)

$\bullet$ Its first singular homology group is torsion free.\\
(See \cite[Corollary 2.2]{Eda2}).

$\bullet$
Its first singular homology group is uncountable;
it contains a divisible torsion-free subgroup of cardinality $2^{\aleph_0}$,
i.e. $\underset{2^{\aleph_0}}{\oplus}\, \mathbb{Q}$.\\
(See \cite[Theorem 3.1]{Eda2} or Corollary~8.7~(a) and Corollary~9.2 of the present paper).

\medskip

Short proofs of some of these statements are presented in \cite[Theorem 2.5]{CC1}.
In~\cite{D}, R.~Distel and P.~Spr{\"u}ssel considered the Freudental compacifications of locally finite connected graphs.
If such a graph has exactly one end, then its Freudental compactification  is homotopy equivalent to the Hawaiian Earring. By analogy to Hawaiian Earring, they proved that the fundamental group
of the Freudental compacification of an arbitrary locally finite connected graph $\Gamma$ can be canonically embedded into the inverse limit of the free groups $\pi_1(\Gamma')$ with $\Gamma'\subseteq \Gamma$ finite.

\bigskip

{\bf Definition 3.3.}
{\it Griffiths' space}
can be constructed as follows: take two Hawaiian Earrings, build
a cone over each of them and finally take their one-point union with respect to the base point of each of the Earrings.



\input epsf
\let\8=\ss
\par\noindent
\hbox{\hskip 4.4truecm{\epsfxsize=8truecm
\epsfbox{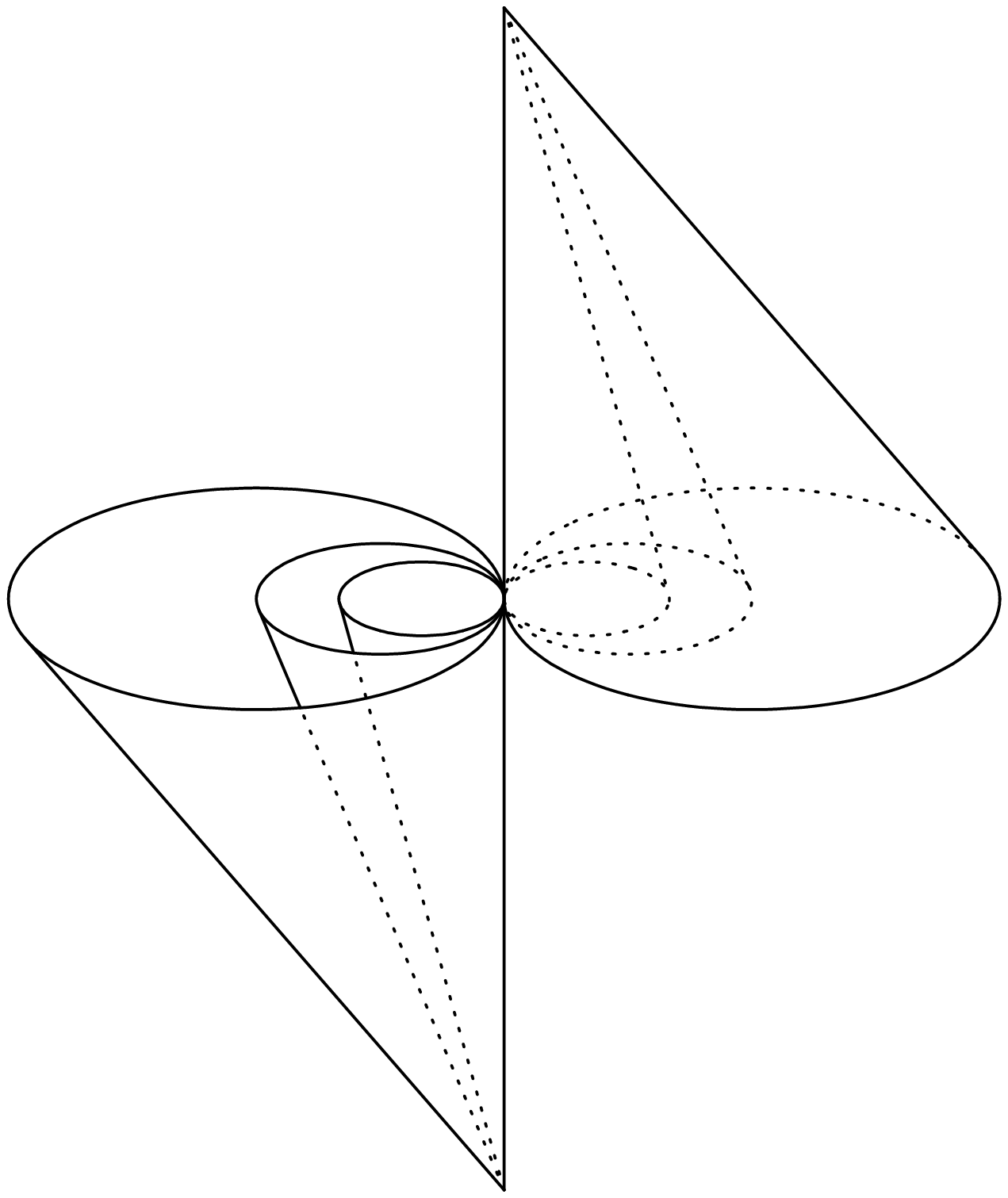}}}
\par\noindent

\begin{center}
Fig.~3
\end{center}

We consider Griffiths' space as embedded into the 3-dimensional Euclidean space so, that both
Hawaiian Earrings lie in the $(x,y)$-plane, have coinciding common tangents and that they are placed in different half-planes with respect to
this tangent. Moreover, the cone-points of
their cones lie in the planes $z=1$ and $z=-1$, respectively (see Fig.~3).
Griffiths' space, if built from two Hawaiian Earrings $H(p_1,p_2,\dots)$ and $H(q_1,q_2,\dots)$, is denoted by $G(p_1,q_1,p_2,q_2,\dots )$, or shortly $G(P;Q)$.
The following theorem is known; see~\cite[formula (3.55)]{Gr1}. For completeness we give here a straightforward proof.

\medskip

{\bf Theorem 3.4.} {\it The fundamental group of Griffiths' space $G(P;Q)$ is isomorphic to the factor group of $\mathcal{W}(P\cup Q)$ by the normal closure of $\langle \mathcal {W}(P),\mathcal{W}(Q) \rangle$.}



\medskip

{\it Proof.} We consider Griffiths' space $Y=G(P;Q)$
with coordinates as in Definition~3.3.
Let $Y_{-1},Y_0,Y_1$ be the subsets of $Y$, consisting of all points with the $z$-coordinate from the intervals  $[-1,-\frac{1}{3}), (-\frac{2}{3}, \frac{2}{3})$, $(\frac{1}{3} , 1]$ respectively.
Note that $\pi_1(Y_{-1})=\pi_1(Y_1)=1$ and $\pi_1(Y_0)\cong \pi_1(H(P\cup Q))$. By the van Kampen Theorem, $\pi_1(Y_{-1}\cup Y_0) $ is the quotient of the group $\pi_1(H(P\cup Q))$ by the normal closure of its subgroup $\pi_1(H(P))$. Finally, $\pi_1(Y)=\pi_1((Y_{-1}\cup Y_0)\cup Y_1)$ is the factor group of the last quotient by the normal closure of the subgroup $\pi_1(H(Q))$. \hfill $\Box$

\medskip

Griffiths' space has the following remarkable properties.

$\bullet$ It is a one-point union of two contractible spaces, but itself it is not contractible.
Moreover its fundamental group is uncountable \cite[Claim 3.4]{Gr1}.

$\bullet$ It is not semi-locally simply connected.

$\bullet$ It does not admit any non-identical covering, and thus it can serve as an example where the univeral covering
space (i.e. the space itself) is not simply connected.\\
(See Corollary 14 in Ch.~2, Sec.~5 of~\cite{Sp} and \cite[Example 2.5.18]{Sp}.)



$\bullet$ It is a space where the natural homomorphism from the fundamental group to the shape group is not an injection \cite[Example 3]{HZ1}.


$\bullet$ Its first singular homology group is torsion free.\\
(This follows from \cite[Theorem 1.2]{Eda3} and \cite[Corollary 2.2]{Eda2}.)

$\bullet$ Its first singular homology group is uncountable; it contains a divisible torsion-free subgroup
of cardinality $2^{\aleph_0}$, i.e. $\underset{2^{\aleph_0}}{\oplus}\, \mathbb{Q}$.\\
(See Corollary~8.7~(b) and Corollary~9.2  of the present paper.)


\medskip

Recall that a space $X$ is called {\it homotopically Hausdorff},
if for every point $x\in X$ and every nontrivial homotopy class
$c\in \pi_1(X,x)$, there is a neighborhood of $x$ which contains no representative for $c$.
In particular, topological manifolds and CW-complexes enjoy this property.
Every one-dimensional space and every planar set is homotopically Hausdorff, see~\cite[Corollary~5.4 (2)]{CC2} and~\cite[Theorem~3.4]{CL}, respectively. In particular, the Hawaiian Earring, is homotopically Hausdorff.
However

$\bullet$ Griffiths' space is not homotopically Hausdorff, since it contains the loop $p_1q_1p_2q_2\dots $,
which is homotopic into any neighborhood of the basepoint, but yet is not nullhomotopic.  This was noticed in \cite[Example~2.0.2 ]{CC2}, but this also follows from our Theorem~8.1.

\medskip

\section{\bf Description of the strategy of the proof}

\medskip


The {\it commutator} of two elements $u,v\in \mathcal{W}(X)$ is the element $u^{-1}\cdot v^{-1}\cdot u\cdot v\in \mathcal{W}(X)$, and it is denoted by $[u,v]$.


\medskip

We describe an idea of the proof of Theorem 8.2 (Theorems 8.1 and 8.3 are proved similarly).
Thus, we consider Griffiths' space $Y=G(P;Q)$.
Let $w_1$, $w_2$ be two reduced tame words over the alphabet $P\cup Q$.
Suppose that
they represent the same element in $H_1(Y)=[\pi_1(Y),\pi_1(Y)]$. Then, by Theorem~3.4 we have
$$w_2= w_1\cdot \prod_{j=1}^m \theta_j^{-1} \cdot \eta_j \cdot \theta_j \cdot \prod_{i=1}^n[\sigma_i,\tau_i]\eqno{(4.1)}$$
for some reduced tame words $\theta_j,\sigma_i,\tau_i,\eta_j\in \mathcal{W}(P\cup Q)$, where for every $j$ holds $\eta_j\in \mathcal{W}(P)$ or $\eta_j\in \mathcal{W}(Q)$.
Therefore $w_2$ is the reduced form of the concatenation
$$w\equiv w_1 \prod_{j=1}^m \theta_j^{-1} \eta_j \theta_j\, \prod_{i=1}^n[\sigma_i,\tau_i]\eqno{(4.2)}$$
By Lemma~2.13, $w$ can be written as
$$w\equiv t_1u_1t_2u_2\dots t_ku_kt_{k+1},\eqno{(4.3)}$$ where the reduced form of each $u_q$ is the empty word and $$w_2\equiv t_1t_2\dots t_{k+1}.\eqno{(4.4)}$$

With the concatenations (4.2)-(4.3) we will associate a geometric object called arch-line-band system (see Section 5). In this system we define leaves and parallelity classes of leaves.
By Proposition~6.6, the number of these classes is finite and we will call it {\it complexity} of concatenations (4.2)-(4.3).

Using transformations (i) and (ii) from Theorem 8.2, we will simplify $w$, $w_1$ and $w_2$ in several steps, so that

1) at each step the new tame words $w_1'$ and $w_2'$ represent the same elements of $H_1(Y)$;

2) the new tame words $w',w_1'$ and $w_2'$ are related by equations, which are analogous to (4.2)-(4.4), have the same parameters $n,m,k$, but smaller complexity:

$$w'\equiv w'_1\prod_{j=1}^m \theta_j'^{-1}\eta'_j\theta'_j\, \prod_{i=1}^n[\sigma'_i,\tau'_i],\eqno{(4.2')}$$
$$w'\equiv t_1'u_1't_2'u_2'\dots t_k'u_k't_{k+1}',\eqno{(4.3')}$$
and
$$w_2'\equiv t_1't_2'\dots t_{k+1}',\eqno{(4.4')}$$
where $t_s',u_l',w_1',\theta_j',\eta_j',\sigma_i',\tau_i'$ are obtained from $t_s,u_l,w_1,\theta_j,\eta_j,\sigma_i,\tau_i$ by removing a finite number of subwords, and
where the reduced form of each $u_q'$ is the empty word;

3) the final word $w_2'$ can be obtained from the final word $w_1'$ by applying a finite number of transformations~(iii).


\section{\bf A visualization of cancellation processes.\\ Arch-line-band systems}

\medskip

In this section we define geometric objects -- arch-line-band systems, which are associated with concatenations
(4.2) and (4.3) of tame words. They are designed to control cancellations in these concatenations.
In Section 7 we introduce transformations of these systems, that will be used in Section 8 to prove our main Theorems 8.1, 8.2 and 8.3.


First we define components of an arch-line-band system,
giving the definition of the entire system on that basis after that.
The following definitions use the visualization
of
countable words as described in Definition 2.2.

\medskip

{\bf (a) Bands.} Here we define two types of bands, line-bands and arch-bands.

\hspace*{5mm}{\it Line-bands.} Let $t$ be a countable word and $T$ be a $t$-segment.
We glue a rectangle $[0,1]\times [0,1]$ to the $t$-segment by a
homeomorphism which identifies the side $[0,1]\times \{0\}$ with this segment. The opposite side $[0,1]\times \{1\}$ is called ${\it free}$.
Inside of this rectangle we draw parallel {\it vertical lines}, starting from the letters of this segment and ending at its free
side. The system, consisting of the $t$-segment, this rectangle and all these lines is called the {\it line-band} associated with $T$. We will say that this line-band is {\it glued} to $T$.

\hspace*{5mm}{\it Arch-bands.} Now let $T$ be a $t$-segment and $T_1$
be a $t^{-1}$-segment. We glue a rectangle $[0,1]\times [0,1]$ to
$T$ and $T_1$ along the sides
$[0,1]\times \{0\}$ and $[0,1]\times \{1\}$ with the help of homeomorphisms which
preserve and reverse the orientation, respectively.
Then we connect the corresponding inverse letters of $T$ and $T_1$ by parallel {\it arcs} inside of this rectangle.
The system, consisting of the segments $T,T_1$, this rectangle and
all these arcs is called the {\it arch-band} associated with $T$ and $T_1$.
We will say that this arch-band {\it connects} $T$ and $T_1$.

\medskip

{\bf (b) Arch-band system.} Let $u$ be a tame word, and suppose that $u$ is reducible to the empty word.
According to Lemma~2.13, there exists a cancelation pattern for $u$,
that is a subdivision $u\equiv z_1z_2\dots z_{2k}$, where some neighboring factors $z_{i(1)},z_{j(1)}$ are mutually inverse,
and after removing them some other neighboring factors $z_{i(2)},z_{j(2)}$ will be mutually inverse, and so on until we get the empty word. The factors $z_{i(s)},z_{j(s)}$ are called {\it associated factors}.

This cancelation procedure may be not unique, but for any fixed one we define the corresponding
arch-band system in the following way.  First we visualize the words $z_i$ on the real line
so that the $z_{i+1}$-segment is situated to the right of the $z_i$-segment in the adjacent position,
$i=1,\dots ,2k-1$. Clearly, the union of these $z_i$-segments is a $u$-segment.

We connect the associated $z_{i(s)}$- and $z_{j(s)}$-segments by an arch-band.
The system, consisting of the $u$-segment and all these arch-bands
is called the {\it arch-band system} associated with $u$.

\medskip


Now, let $w$ be a tame word, which is represented as a concatenation in two ways, see (4.2) and (4.3),
where $u_i$'s are tame words, which are reducible to the empty word.
We visualize $w$ by a $w$-segment.

\medskip

{\bf (c) Upper arch-line-band system.}
Consider $2k+1$ consecutive subsegments of the $w$-segment, which correspond
to the factors on the right hand side of (4.3).
For every $t_i$-segment we draw the associated line-band and for every $u_i$-segment we draw the associated arch-band system.
The system, consisting of the $w$-segment, the line-bands and the arch-bands is called the {\it upper arch-line-band system} associated with the concatenation (4.3). The corresponding bands are called {\it upper}.

\medskip



{\bf (d) Lower arch-line-band system.}
Consider $1+3m+4n$ consecutive subsegments of the $w$-segment,
which correspond to the factors in the right side of (4.2).
Then for each $i,j$ we

(i) connect the $\theta_j^{-1}$- and  $\theta_j$-segments by an arch-band,

(ii) connect the $\sigma_i^{-1}$- and $\sigma_i$-segments
and also connect the $\tau_i^{-1}$- and  $\tau_i$-segments by arch-bands,

(iii) glue line-bands to the $w_1$- and to the $\eta_j$-segments.

\medskip

The system, consisting of the $w$-segment, all these line-bands and arch-bands is called the
{\it lower arch-line-band system} associated with the concatenation (4.2). The corresponding bands are called {\it lower}.

\medskip

{\bf (e) The arch-line-band system} associated with the equations (4.2) and (4.3) is the union of the corresponding upper and lower arch-line-band systems based on the same $w$-segment (see Fig.~4). It will be denoted by

$$\frac{t_1u_1t_2u_2\dots t_ku_kt_{k+1}}{w_1\prod_{j=1}^m\theta_j^{-1}\eta_j \theta_j \prod_{i=1}^n[\sigma_i,\tau_i]}.\eqno{(5.1)}$$

\medskip

\input{ArchLine5.pic}
\vspace*{-12mm}
\begin{center}
Fig.~4
\end{center}

{\bf Definition.} Let $(\dots, \sigma_{i-1},\sigma_i,\sigma_{i+1},\dots )$ be a finite or infinite sequence of arches or vertical lines of our system, which lie alternately in upper or lower bands such that
the terminal point of each $\sigma_j$ coincides with the initial point of its successor $\sigma_{j+1}$.
The union of these vertical lines and arches is called a {\it curve}.

\medskip

\input{Leaf6.pic}

\begin{center}
Fig.~5
\end{center}

A {\it leaf} in the arch-line-band system is a maximal curve (see Fig.~5).
Clearly, vertical lines can appear only at the beginning or at the end of the leaf.
Note that a leaf cannot intersect itself and it never splits.

\section{\bf Properties of leaves of an arch-line-band system}

\medskip



 {\bf Proposition 6.1.} {\it Any leaf of an arch-line-band system consists of a finite number of arcs and vertical lines.}


 {\it Proof.} 
We will use notations of Section 5.
 Any leaf intersects the $w$-segment of the system at points, which are marked by
 letters of $\{x,x^{-1}\}$ for some $x\in X$.
 Since $w$ is a tame word, each letter of $X^{\pm}$ occurs only finitely many times in $w$.
 Thus, any leaf intersects the $w$-segment in a finite set of points. \hfill $\Box$

\medskip

Therefore each leaf $\sigma$ can be written in the form $\sigma=\sigma_1\sigma_2\dots\sigma_p$, where $p$ is finite and all $\sigma_i$ are arches, except for $\sigma_1$ and $\sigma_p$, which also may be vertical lines.

\medskip

Now we introduce the concept of parallelity for arches, lines and leaves within arch-line-band systems. This concept  technically differs from the corresponding one for arches within arch-systems (cf. Definition~2.3.1)--2)), though it is morally in the same spirit.

\medskip

{\bf Definition 6.2.}  a) Two arches (lines)
of an arch-line-band system
are called {\it parallel} if they lie in the same 
arch-band (line-band).

b) Two nonclosed leaves $\sigma_1\sigma_2\dots\sigma_p$ and $\tau_1\tau_2\dots\tau_q$ are called {\it parallel} if $p=q$ and $\sigma_i,\tau_i$ are parallel for any $i=1,\dots, p$.

Two closed leaves $\sigma_1\sigma_2\dots\sigma_p$ and $\tau_1\tau_2\dots\tau_q$ are called {\it parallel} if $p=q$ and there is a natural $k$, such that the arches $\sigma_i,\tau_{i+k}$ are parallel for any $i=1,\dots, p$ (addition modulo $p$).

\medskip

{\bf Remark 6.3.} Let $\sigma_1\sigma_2\dots\sigma_p$ and $\tau_1\tau_2\dots\tau_q$ be two nonclosed leaves such that
$p\leqslant q$ and $\sigma_i,\tau_i$ are parallel for any $i=1,\dots ,p$. Then $\sigma_p$ ends on the free side of some
line-band and so does $\tau_p$. Therefore $p=q$.

\medskip

{\bf Definition 6.4.} Let $S$ be a topological space with a distinguished subset $D$. Two non-closed curves $l_1,l_2$ in $S$ with endpoints in $D$ will be called {\it homotopic relatively $D$} if there exists a homotopy carrying $l_1$ to $l_2$, which moves the endpoints of $l_1$ and $l_2$ inside of $D$.
If we apply this term to two closed curves, we mean that they are freely homotopic.

\medskip

Sometimes we will consider an arch-line-band system as a topological space, which is the union of
the corresponding rectangles glued to a segment of the real line.

\medskip

{\bf Lemma 6.5.} {\it Let $A$ be an arch-line-band system and $D=b_1\cup \dots \cup b_t$ be the union of the free sides of its line-bands. Two leaves of $A$ are parallel if and only if they are homotopic relatively~$D$.}

\medskip

{\it Proof.} If two leaves of $A$ are parallel, then, clearly, they are homotopic relatively $D$.
We prove the converse statement.

Let $R$ be a rose with one vertex $v$, petals corresponding to the arch-bands and with thorns corresponding to the line-bands of $A$ (see Fig.~6).

\medskip

\input{SurfaceThorn7.pic}

\begin{center}
Fig.~6
\end{center}

There is a natural homotopy equivalence $\pi:A\rightarrow R$, which sends the $w$-segment to $v$, arch-bands to petals and line-bands to thorns. Free sides of line-bands go to the vertices of thorns different from $v$. The set of these vertices is $\pi(D)$.

Let $l_1,l_2$ be two nonparallel leaves in $A$, both closed or both nonclosed. If they are closed, their images in $R$ are different cyclic cyclically reduced paths. If $l_1,l_2$ are nonclosed, their images in $R$  are different reduced paths with endpoints in $\pi(D)$. In both cases these paths are not homotopic in $R$ relatively $\pi(D)$, and so the leaves are not homotopic in $A$ relatively $D$. \hfill $\Box$

\bigskip

{\bf Proposition 6.6.} {\it
There is only a finite number of parallelity classes of leaves in the arch-line-band system $(5.1)$}.

\medskip

{\it Proof.} Let $A$ be the topological space associated with the system (5.1).
We embed this space into a surface, so that this embedding will be a homotopy equivalence. To construct the surface we will glue to $A$ a finite number of triangles and rectangles.

More exactly, for any two adjacent upper (lower) bands in $A$ their non-free sides are either disjoint or meet at a point of the $w$-segment; in the first case we glue a rectangle, in the second case we glue a triangle as it shown in Fig.~7. Denote the resulting space by $T$
and note that it is a  finite genus orientable surface.

\medskip

\input{Surface8.pic}

\begin{center}
Fig.~7
\end{center}

The boundary of $T$ contains parts of the lateral sides of all bands, some sides of the triangles and rectangles, and
the free sides of all line-bands.
Within $\partial T$ we consider these free sides $b_1,\dots ,b_t $ as {\it distinguished boundary segments}.

Since the space $A$ is contained in $T$, we may assume that all leaves of $A$ are contained in $T$. These leaves are of two types: closed curves in the interior of $T$ and curves with ends on distinguished boundary segments.

Taking into account that $A$ and $T$ are homotopy equivalent relatively $D=b_1\cup \dots \cup b_t$,
we conclude from Lemma 6.5 that two leaves of $A$ are parallel if and only if they are homotopic in $T$ relatively $D$. Now the statement follows from Lemma 6.7. \hfill $\Box$

\medskip

For any compact surface $S$ we denote by $g(S)$ its genus and by $b(S)$ the number of its boundary components.

\medskip

{\bf Lemma 6.7.} {\it Let $S$ be a compact surface with a finite set of distinguished segments $b_1,\dots ,b_t$
on its boundary. Let $L$ be a set of disjoint simple curves in $S$, such that each curve of $L$ is either

{\rm (1)} a closed curve in the interior of $S$, or

{\rm (2)} a nonclosed curve, which has both endpoints in $\{b_1\cup \dots \cup b_t\}$ and lies,
apart from them, in the interior of $S$.


Then the number of relatively-$\{b_1\cup \dots \cup b_t\}$ homotopy classes of these curves does not exceed a constant depending only on $g(S)$,  $b(S)$ and $t$.}

\medskip

{\it Proof.} We will proceed by induction on $(g(S), b(S))$ trying to decrease the pair in the lexicographical order  $\preccurlyeq$ on $\mathbb{N}\times \mathbb{N}$ by cutting $S$ along a curve $l\in L$.
There are 2 cases.

{\it Case 1.} $l$ is a nonseparating curve. Then, after cutting along $l$, we get a new surface $S_1$. In this case there are 3 subcases.

\hspace*{5mm} a) $l$ is a closed curve. Then  $g(S_1)=g(S)-1$ and $b(S_1)=b(S)+2$.

\hspace*{4.5mm} b) $l$ is a nonclosed curve with endpoints at different boundary components of $S$. Then
$g(S_1)=g(S)$ and $b(S_1)=b(S)-1$.

\hspace*{4.5mm} c) $l$ is a nonclosed curve with endpoints at the same boundary component of $S$. Then
$g(S_1)=g(S)-1$ and $b(S_1)=b(S)+1$.

{\it Case 2.} $l$ is a separating curve. Then, if we cut along $l$, we get two surfaces $S_1,S_2$. In this case there are 2 subcases.

\hspace*{5.4mm} a) $l$ is a closed curve. Then $g(S_1)+g(S_2)=g(S)$ and $b(S_1)+b(S_2)=b(S)+2$.

\hspace*{5mm} b) $l$ is a nonclosed curve with endpoints at the boundary of $S$. Since $l$ is se\-pa\-rating, its  endpoints lie at the same boundary component of $S$. Then $g(S_1)+g(S_2)=g(S)$ and $b(S_1)+b(S_2)=b(S)+1$, and $b(S_1)>0$, $b(S_2)>0$.

\medskip

We see, that in Case 1 the pair $(g(\cdot ),b(\cdot ))$ decreases in the lexicographical order $\preccurlyeq$.
In Case~2 we have also $(g(S_j),b(S_j)) \prec (g(S),b(S))$ for $j=1,2$ except of three subcases, where for some $j$ holds

\medskip

(i) $l$ is closed and $(g(S_j),b(S_j))$ equals to $(0,1)$,

(ii) $l$ is closed and $(g(S_j),b(S_j))$ equals to $(0,2)$,

(iii) $l$ is nonclosed and $(g(S_j),b(S_j))$ equals to $(0,1)$.

\medskip

Note, that the curve $l$ lies on $\partial S_j$.

In Subcase (i), $S_j$ is a disk and $l=\partial S_j$. In particular $l$ is null-homotopic in $S$.

In Subcase (ii), $S_j$ is an annulus and $l$ is one of its boundary components. In particular $l$ is parallel to the other boundary component of $S_j$, which is a boundary component of $S$.

In Subcase (iii), $S_j$ is a disc and $l$, having endpoints on distinguished segments of $\partial S$, is parallel to $\partial S$.

\medskip

Thus, if $L$ contains more than $b(S)+2t^2+1$ relatively-$\{b_1\cup \dots \cup b_t\}$ homotopy classes of curves, then one of them does not fall into the  exceptional subcases (i), (ii) or (iii). Cutting $S$ along this curve, we get one or two surfaces with smaller $(g(\cdot),b(\cdot))$. Note that for each of the resulting surfaces $g(\cdot )$ cannot increase, $b(\cdot)$ can increase by at most 2 and $t$ can increase by at most 2. Thus for each of the new surfaces we can apply the induction.

Note, that the number of the relatively-$\{b_1\cup \dots \cup b_t\}$ homotopy classes of curves in $L$ does not exceed the sum of the corresponding numbers after cutting of $S$ along a curve $l\in L$.
This completes the proof.
\hfill $\Box$


\section{\bf Transformations of arch-line-band systems}

Now we consider an important transformation of an arch-line-band system:

{\bf removing the parallelity class of a leaf}. Let $\sigma$
be a leaf in the arch-line-band system (5.1)
and $[\sigma]$ be the parallelity class of $\sigma$.
We then ``remove $[\sigma ]$'', i.e.\ we remove from the factors of (5.1) all letters that lie on the intersection of the $w$-segment and leaves from $[\sigma ]$, and we also remove these leaves from our picture. Due to Proposition~6.6, only finitely many subwords will be removed from each of the factors $t_s,u_l,w_1,\theta_j,\eta_j,\sigma_i,\tau_i$ of (5.1). The remaining parts of these factors form new factors  $t_s',u_l',w_1',\theta_j',\eta_j',\sigma_i',\tau_i'$ and we obtain a new system:
$$\frac{t_1'u_1't_2'u_2'\dots t_k'u_k't_{k+1}'}{w_1'\prod_{j=1}^m\theta_j'^{-1}\eta_j'\theta_j' \prod_{i=1}^n[\sigma_i',\tau_i']}.\eqno{(7.1)}$$

Note, that the remaining letters are still connected by the same upper or lower arches as before.
Therefore we conclude:

\medskip

{\bf Remark 7.1.} 1) The words $u'_s$ are  reducible to the empty word.

2) We have ${\theta^{-1}_j}' \equiv {\theta'_j}^{-1}$, ${\sigma^{-1}_i}' \equiv {\sigma'_i}^{-1}$,
${\tau^{-1}_i}' \equiv {\tau'_i}^{-1}$, i.e.\ letters were deleted from pairs
of inverse factors in such a way, that the new factors remained inverse to each other.\par\noindent


3) By Remark 2.12, all the words in the numerator and the denominator of (7.1) are tame, since they are finite concatenations of subwords of tame words.

4) Summing up 1)--3) we obtain: the system (7.1) indeed satisfies
all conditions of an arch-line-band system, and has complete analogous structure of factors
as (5.1) (although some factors might have become empty). In particular, the numerator and the denominator in (7.1)
are equal.

\medskip

We are especially interested in the exact description of the structure of $t_1',\dots ,t_k'$ and $w_1'$. This structure depends on the type of the parallelity class of $\sigma$.


{\bf (1)} If $\sigma$ is closed it cannot run through $w_1$ and $t_i$. In this case we have $w_1\equiv w_1'$ and $t_i\equiv t_i'$ for every $i$.

\medskip

{\bf (2)} If $\sigma$ starts at an upper $t_i$-band and ends at another upper $t_j$-band, then there exists a tame word $W$, such that $t_i$ contains $W$, $t_j$ contains $W^{-1}$ and we have the following transformation:

$$\begin{array}{ll} t_i\equiv A_iWB_i & \mapsto \hspace*{2mm} t_i'\equiv A_iB_i,\vspace*{2mm}\\
t_j\equiv A_jW^{-1}B_j &\mapsto \hspace*{2mm} t_j'\equiv A_jB_j,\vspace{2mm}\\
t_l & \mapsto \hspace*{2mm} t_l'\equiv t_l \hspace*{10mm} {\text{\rm for}}\quad l\neq i,j,\vspace*{2mm}\\
w_1 & \mapsto \hspace*{2mm} w_1'\equiv w_1.
\end{array}$$

\medskip

{\bf (2')} If $\sigma$ starts and ends at the same upper $t_i$-band, then $t_i$ has the form $t_i\equiv A_iWB_iW^{-1}C_i$ and we have the following transformation:

$$\begin{array}{ll} t_i\equiv A_iWB_iW^{-1}C_i & \mapsto \hspace*{2mm} t_i'\equiv A_iB_iC_i,\vspace*{2mm}\\
t_l & \mapsto \hspace*{2mm} t_l'\equiv t_l \hspace*{10mm} {\text{\rm for}}\quad l\neq i,\vspace*{2mm}\\
w_1 & \mapsto \hspace*{2mm} w_1'\equiv w_1.
\end{array}$$

\medskip

{\bf (3)} If $\sigma$ starts at a lower $\eta_i$-band and ends at an upper $t_j$-band, then there is a common subword $W$ of $\eta_i$ and $t_j$, such that we have the following transformation:

$$\begin{array}{ll} t_j\equiv A_jWB_j & \mapsto \hspace*{2mm} t_j'\equiv A_jB_j,\vspace*{2mm}\\
t_l & \mapsto \hspace*{2mm} t_l'\equiv t_l \hspace*{10mm} {\text{\rm for}}\quad l\neq j,\vspace*{2mm}\\
w_1 & \mapsto \hspace*{2mm} w_1'\equiv w_1.
\end{array}$$

\medskip

{\bf (4)} If $\sigma$ starts at the lower $w_1$-band and ends at an upper $t_j$-band, then there is a common subword $W$ of $w_1$ and $t_j$, such that we have the following transformation:

$$\begin{array}{ll} t_j\equiv A_jWB_j & \mapsto \hspace*{2mm} t_j'\equiv A_jB_j,\vspace*{2mm}\\
t_l & \mapsto \hspace*{2mm} t_l'\equiv t_l \hspace*{10mm} {\text{\rm for}}\quad l\neq j,\vspace*{2mm}\\
w_1\equiv AWB & \mapsto \hspace*{2mm} w_1'\equiv AB.
\end{array}$$

\medskip

{\bf (5)} If $\sigma$ starts at a lower $\eta_i$-band and ends at the lower $w_1$-band, then there exists a tame word $W$, such that $\eta_i$ contains $W$, $w_1$ contains $W^{-1}$ and we have the following transformation:

$$\begin{array}{ll}
t_l & \mapsto \hspace*{2mm} t_l'\equiv t_l \hspace*{10mm} {\text{\rm for}}\quad l=1,\dots,k+1,\vspace*{2mm}\\
w_1\equiv AW^{-1}B & \mapsto \hspace*{2mm} w_1'\equiv AB.
\end{array}$$

\medskip

{\bf (6)} If $\sigma$ starts at a lower $\eta_i$-band and ends at a lower $\eta_j$-band, then $w_1'\equiv w_1$
and $t_l'\equiv t_l$ for all~$l$.

\medskip

{\bf (7)} If $\sigma$ starts and ends at the lower $w_1$-band, then $w_1$ has the form $w_1\equiv AWBW^{-1}C$ and we have the following transformation:

$$\begin{array}{ll}
t_l & \mapsto \hspace*{2mm} t_l'\equiv t_l \hspace*{10mm} {\text{\rm for}}\quad l=1,\dots,k+1,\vspace*{2mm}\\
w_1\equiv AWBW^{-1}C & \mapsto \hspace*{2mm} w_1'\equiv ABC.
\end{array}$$

\medskip

\section{\bf Main theorems}


\medskip

{\bf Theorem 8.1.} {\it Let $Y$ be Griffiths' space $G(p_1,q_1,p_2,q_2,\dots )$. Two reduced tame words $w_1,w_2$ over the alphabet $\{ p_1,q_1,p_2,q_2,\dots \}$ represent the same element in $\pi_1(Y)$ if and only if each of them can be transformed to the same reduced tame word by applying a finite number of transformations of type {\rm (i)} and finally the transformation of type~{\rm (ii)}:}

\medskip

{\rm (i)}  {\it Deletion of a countable subword, containing only $p$- or $q$-letters;}

{\rm (ii)} {\it Reducing the resulting countable word.}

\bigskip

{\it Proof.}  The ``if''-part follows from Theorem~3.4: by this theorem any restricted countable word containing only $p$- or $q$-letters represents the identity element of $\pi_1(Y)$. So, we prove the ``only~if''-part.

Let $w_1,w_2$ be two reduced tame words over the alphabet $\{ p_1,q_1,p_2,q_2,\dots \}$
representing the same element in $\pi_1(Y)$.
By  Theorem~3.4, the word $w_2$ is the reduced form of a finite concatenation
$$w_1\prod_{j=1}^m \theta_j^{-1}\eta_j \theta_j,\eqno{(8.1)}$$
where $\theta_j\in \mathcal{W}(p_1,q_1,p_2,q_2,\dots )$
and each $\eta_j\in \mathcal{W}(p_1,p_2,\dots)$ or $\eta_j\in
\mathcal{W}(q_1,q_2,\dots)$.
By Lemma~2.13, there exists
a decomposition of (8.1),
$$w_1\prod_{j=1}^m\theta_j^{-1}\eta_j^{-1}\theta_j\equiv t_1u_1t_2u_2\dots t_ku_kt_{k+1},\eqno{(8.2)}$$
where all the subwords $u_i$ are reducible to the empty word, all $t_i$ are irreducible and
$$w_2\equiv t_1t_2\dots t_{k+1}.$$
Now we consider the arch-line-band system
$$\frac{t_1u_1t_2u_2\dots t_ku_kt_{k+1}}{w_1\prod_{j=1}^m\theta_j^{-1}\eta_j\theta_j}.\eqno{(8.3)}$$
Consecutively remove all parallelity classes of leaves from it, which start at lower $\eta_j$-bands (cf. Section~7).
This process is finite by Proposition~6.6.
Then the words $\eta_j$ completely disappear and we get a new arch-line-band system
$$\frac{t_1'u_1't_2'u_2'\dots t_k'u_k't_{k+1}'}{w_1'\prod_{j=1}^m\theta_j'^{-1}\theta_j'}.\eqno{(8.4)}$$
By analogy with $w_2\equiv t_1t_2\dots t_{k+1}$ we denote $w_2'\equiv t_1't_2'\dots t'_{k+1}$. According to Cases (3), (5), (6) from Section 7, the countable words $w_2'$ and $w_1'$ are obtained from $w_2$ and $w_1$ by applying a finite number of transformations (i).

Since the numerator and the denominator in (8.4) are equal as countable words and still tame (see statements 4) and 3) of Remark 7.1), they have the same reduced form by Theorem~2.9. Hence $w_2'$ and $w_1'$
have the same reduced form. In other words, after applying to $w_2'$ and $w_1'$ the transformation~(ii), we get the same reduced tame word. \hfill $\Box$

\bigskip

{\bf Theorem 8.2.} {\it Let $Y$ be Griffiths' space $G(p_1,q_1,p_2,q_2,\dots )$. Two reduced tame words $w_1,w_2$ over the alphabet $\{ p_1,q_1,p_2,q_2,\dots \}$ represent the same element in $H_1(Y)$ if and only if each of them can be transformed to the same tame word
by applying a finite number of transformations of type {\rm (i)}, then of  type {\rm (ii)} and finally of type {\rm (iii)}:}

\medskip

(i)  {\it Deletion of a countable subword, containing only $p$- or $q$-letters;}

(ii) {\it Deletion of two distinct mutually inverse countable subwords: $AXBX^{-1}C\rightarrow ABC$;}

(iii) {\it Permuting two consecutive countable subwords}.

\medskip

{\it Proof.} We will prove the nontrivial ``only if'' part. Let $w_1,w_2$ be two reduced tame words over the alphabet $\{ p_1,q_1,p_2,q_2,\dots \}$ representing the same element in $H_1(Y)$.
By  Proposition~3.4, the word $w_2$ is the reduced form of a finite concatenation
$$w\equiv w_1 \prod_{j=1}^m \theta_j^{-1}\eta_j\theta_j\, \prod_{i=1}^n[\sigma_i,\tau_i],\eqno{(8.5)}$$
where each $\eta_j$ contains only $p$- or $q$-letters.
This means, that there exists a decomposition of (8.5),
$$w_1 \prod_{j=1}^m \theta_j^{-1}\eta_j\theta_j\, \prod_{i=1}^n[\sigma_i,\tau_i]\equiv t_1u_1t_2u_2\dots t_ku_kt_{k+1},\eqno{(8.6)}$$
where all the subwords $u_i$ are reducible to the empty word and

$$w_2\equiv t_1t_2\dots t_{k+1}.$$

Now we consider the arch-line-band system
$$\frac{t_1u_1t_2u_2\dots t_ku_kt_{k+1}}{w_1\prod_{j=1}^m\theta_j^{-1}\eta_j\theta_j\, \prod_{i=1}^n[\sigma_i,\tau_i]}.\eqno{(8.7)}$$
We will simplify this system by removing certain parallelity classes of leaves.
First we eliminate all $\eta_j$ exactly in the same way as in the proof of Theorem~8.1 (so, we use transformations (i) several times). Thus, without loss of generality, we may assume that we have the arch-line-band system

$$\frac{t_1u_1t_2u_2\dots t_ku_kt_{k+1}}{w_1\prod_{j=1}^m\theta_j^{-1}\theta_j\, \prod_{i=1}^n[\sigma_i,\tau_i]}.\eqno{(8.8)}$$

The only  line-bands of this system are the upper $t_s$-bands, $s=1,\dots ,k+1$, and the lower $w_1$-band.
Note that $t_1t_2\dots t_{k+1}$ is a tame word, but may be already non-reduced.

Now suppose that the system contains a leaf $\sigma$ starting at some $t_p$-band
and ending at the same or another $t_q$-band.
By Cases (2) and (2') from Section~7, the tame word $w_2\equiv t_1t_2\dots t_{k+1}$ has the form $AXBX^{-1}C$.
Removing the parallelity class of $\sigma$, we do not change the form of~$(8.8)$ and do not change $w_1$,
but we replace $w_2$ by $ABC$ (this is the transformation (ii)).


We remove all parallelity classes starting and ending at the union of the upper $t$-bands. Analogously we remove all parallelity classes starting and ending at the lower $w_1$-band.

Thus we may assume, that every parallelity class starting at an upper $t$-band ends at the lower $w_1$-band. Conversely, every parallelity class starting at the lower $w_1$-band ends at some upper $t$-band. Moreover, every letter of $w_2\equiv t_1t_2\dots t_{k+1}$ and every letter of $w_1$ lies on some leaf.
Therefore there is a partition $w_2\equiv v_1v_2\dots v_r$ and a permutation $\sigma\in S_r$,  such that $w_1\equiv v_{1\sigma}v_{2\sigma}\dots v_{r\sigma}$. Thus, $w_1$ and $w_2$ can be made equal by a finite number of transformations (iii). \hfill $\Box$
\bigskip

{\bf Theorem 8.3.} {\it Let $Z$ be the Hawaiian Earring $H(p_1,p_2,\dots )$. Two reduced tame words over the alphabet $\{ p_1,p_2,\dots \}$ represent the same element in the first homology group $H_1(Z)$ if and only if each of them can be transformed to the same tame word
by applying a finite number of transformations of type {\rm (ii)} and finally of type {\rm (iii)}:}

\medskip

(ii) {\it Deletion of two distinct mutually inverse countable subwords: $AXBX^{-1}C\rightarrow ABC$;}

(iii) {\it Permuting two consecutive countable subwords}.

\medskip

{\it Proof.}
We will prove the nontrivial ``only if'' part.
The word $w_2$ is the reduced form of a finite concatenation
$$w_1\, \prod_{i=1}^n[\sigma_i,\tau_i]\eqno{(8.9)}$$ (compare with (8.5)).
The further proof is the same as the proof of Theorem 8.2, starting from formula (8.8) with empty $\theta_j$. \hfill $\Box$

\medskip

Note that this theorem can be also deduced from~\cite[Theorem~6.1]{CC2} and conversely.



\medskip

We fix two disjoint infinite countable alphabets $P=\{p_1,p_2\dots \}$ and $Q=\{q_1,q_2,\dots \}$
for the rest of this section. Below we introduce legal and $(P,Q)$-legal words, which, in view of Lemma~8.6,
will be used later to construct uncountably many words with certain nice properties.

\medskip


{\bf Definition 8.4.}

1) A countable word $U$ over the alphabet $P$ is called {\it legal}, if it is restricted and for any presentation $U\equiv AXBX^{-1}C$ the subword $X$ is finite.


2) A countable word $U$ over the alphabet $P\cup Q$ is called {\it $(P,Q)$-legal} if
it is restricted,
every subword of $U$ containing only $p$- or only $q$-letters is finite,
and for any presentation $U\equiv AXBX^{-1}C$ the subword $X$ is finite.


\medskip

The proof of the following useful lemma is straightforward.

{\bf Lemma 8.5.} {\it The class of legal words over the alphabet $P$ and the class of $(P,Q)$-legal words over the
alphabet $P\cup Q$ both are closed under applying the transformations {\rm (i)} and {\rm (ii)} from {\rm Theorem~8.2}
a finite number of times.}











\medskip

{\bf Lemma 8.6.} {\it
{\rm 1)} A reduced legal word $U$ over the alphabet $P$ determines a nonzero element in $H_1(Z)$ if and only if it is infinite or it is finite and has a nonzero exponent sum in some letter.

{\rm 2)} A reduced $(P,Q)$-legal word $U$ over the alphabet $P\cup Q$ determines a nonzero element in $H_1(Y)$ if and only if it is infinite.}


\medskip

{\it Proof.} We will prove only the more difficult statement 2).
By Theorem 8.2, $U$ determines the zero element in $H_1(Y)$ if and only if it can be carried
to the empty word by applying a finite number of transformations (i), (ii). By Lemma~8.5, each transformation leaves words in the class of $(P,Q)$-legal words. Therefore at each step at most two finite subwords can be deleted.
Thus, if $U$ is infinite, we cannot get the empty word in a finite number of steps. If $U$ is finite, we can delete
all its letters with the help of several transformations~(i).

The proof of the statement 1) is analogous and uses Theorem~8.3. \hfill $\Box$


\medskip

Recall, that $\aleph_0$ denotes the cardinality of natural numbers, and so $2^{\aleph_0}$ is the cardinality of continuum.

\medskip

{\bf Corollary 8.7.} {\it For the Hawaiian Earring $Z=H(p_1,p_2,\dots )$ and for Griffiths' space $Y=G(p_1,q_1,p_2,q_2,\dots )$, there are
$2^{\aleph_0}$
functions $k:\mathbb{N}\rightarrow \mathbb{N}$, such that

{\rm (a)} all the elements $[p_{k(1)},p_{k(2)}][p_{k(3)},p_{k(4)}]\dots $ are different in $H_1(Z)$;

\vspace*{1mm}

{\rm (b)} all the elements $[p_{k(1)},q_{k(1)}][p_{k(2)},q_{k(2)}]\dots $ are different in $H_1(Y)$.}

\medskip

{\it Proof.} Consider the set of all strictly monotone functions $k:\mathbb{N}\rightarrow \mathbb{N}$. It has cardinality $2^{\aleph_0}$.
Two such functions $f:\mathbb{N}\rightarrow \mathbb{N}$ and $g:\mathbb{N}\rightarrow \mathbb{N}$
will be called {\it equivalent}, if there exist natural $n,m$ such that for
every $l\geqslant 0$ we have that
$f(n+l)=g(m+l)$. Roughly speaking, the tails of equivalent
functions are the same.
Since each equivalence class of these functions is countable,
the cardinality of the set of equivalence classes is $2^{\aleph_0}$. We choose a system $\mathcal{K}$ of representatives of these classes to define $2^{\aleph_0}$ infinite products according to (a) and (b), respectively.

It is easy to see that for any two such products $U,V$, the concatenation $UV^{-1}$ is reduced legal in Case (a) or
reduced $(P,Q)$-legal in Case (b).
The proof can be completed by applying Lemma~8.6. \hfill $\Box$

\section{Divisible subgroups of  $H_1(Z)$, $H_1(Y)$ and $\pi_1(Y)$}

\medskip

\hspace*{6mm}{\bf Definition 9.1.} An element $g$ of a group $G$ is called
{\it {$\mathbb{N}$}-divisible} if for any $n\in \mathbb{N}$,
there exists an element $x\in G$ such that $g=x^n$.  An abelian group $G$ is called {\it divisible} if every element of $g\in G$ is {$\mathbb{N}$}-divisible.

\medskip

{\bf Corollary 9.2.} {\it The homology groups $H_1(Z)$ and $H_1(Y)$ of the Hawaiian Earring $Z=H(p_1,p_2,\dots )$ and of Griffiths' space $Y=G(p_1,q_1,p_2,q_2,\dots )$ contain subgroups isomorphic to $\underset{2^{\aleph_0}}{\oplus}\, \mathbb{Q}$.

Furhermore, the fundamental group of Griffith's space contains a subgroup isomorphic to $\mathbb{Q}$.}


\medskip


{\it Proof.} a) We construct a nontrivial $\mathbb{N}$-divisible element in $H_1(Y)$ using limits of word sequences.
Let $a_1,a_2,\dots $ be an infinite sequence
of countable words. First we define a countable word $w$ as a limit word\footnote{This limit word is a modification of those described in \cite[pages 235-236]{CC1} and \cite[page 14]{Z3}.} resulting from the following recursive process: set $w_1=a_1$ and construct
$w_{i+1}$ from $w_i$ by inserting $a_{i+1}^{i+1}$ after each
occurrence of $a_i$:

$$\begin{array}{lllllll}
w_1= & a_1 &  &  &  &  &\\
w_2= & a_1 & a_2 &  & a_2 &  & \\
w_3= & a_1 & a_2 & a_3a_3 a_3 & a_2 & a_3a_3a_3 & \\
\vdots &&&&&&
\end{array}
$$


\medskip

Specifically for $a_i=[p_i,q_i]$, $i=1,2,\dots$, we denote the limit word $w$ by $W$. Clearly, $W$ is restricted and reduced, and so $W$ is a nontrivial element of $G=\mathcal{W}(P\cup Q)$. By Theorem~3.4,
there is a natural homomorphism $G\rightarrow \pi_1(Y)$.
We compose it with the canonical homomorphism $\pi_1(Y)\rightarrow H_1(Y)$
and denote the image of $W$ in $H_1(Y)$ by $\overline{W}$.

Note, that if we delete in $W$ a finite number of $a_i$'s, the resulting word will not change modulo $[G,G]$.
Thus, by deleting all occurrences of $a_1$ in $W$, we see that $W$ is a square modulo $[G,G]$. Similarly, for any natural $n$, if we delete all occurrences
of $a_1,a_2,\dots ,a_n$ in $W$, we see that $W$ is an $(n+1)$-th power
in $G$ modulo $[G,G]$.
Therefore $\overline{W}$ is $\mathbb{N}$-divisible in $H_1(Y)$. Moreover, $\overline{W}$ is nontrivial in $H_1(Y)$ by Lemma 8.6.2).

\medskip

b) Now we construct
$2^{\aleph_0}$ $\mathbb{N}$-divisible elements in $H_1(Y)$. For that we will use functions
$k:\mathbb{N}\rightarrow \mathbb{N}$ from the class  $\mathcal{K}$, which was defined in the proof of Corollary~8.7.
For each function $k\in \mathcal{K}$, we define the limit word $W_k$ just by putting
$a_i=[p_{k(i)},q_{k(i)}]$ in the above construction. And as above, one can show, that these words determine $\mathbb{N}$-divisible elements $\overline{W_k}$ in $H_1(Y)$.

It remains to prove, that the words $W_k$, $k\in \mathcal{K}$, determine
different elements in $H_1(Y)$. This follows from Lemma~8.6.2), if we prove that for different $k,l\in \mathcal{K}$ the reduced form $R_{k,l}$
of the word $W_k^{-1}W_l$ is infinite and $(P,Q)$-legal. The word $R_{k,l}$ is infinite, since
$$W_k^{-1}W_l\equiv (\dots [p_{k(3)},q_{k(3)}]^{-1}[p_{k(2)},q_{k(2)}]^{-1}[p_{k(1)},q_{k(1)}]^{-1}
[p_{l(1)},q_{l(1)}][p_{l(2)},q_{l(2)}][p_{l(3)},q_{l(3)}]\dots),$$
and so the reduction can affect only a finite number of the middle commutators.

The word $R_{k,l}$ is $(P,Q)$-legal, since

(1) $p$- and $q$-subwords of $R_{kl}$ are finite;

(2) if $R_{k,l}\equiv AXBX^{-1}C$, then $X$ is finite.

The first assertion is clear, since $p$- and $q$-subwords of $W_k$ and $W_l$ contain only one letter.
For the second one we need the following claim, which can be proved straightforward.

\medskip

{\sl Claim.} {\it Every infinite subword of the limit word $w$ contains an infinite subword of the form
$a_na_{n+1}a_{n+2}\dots $, which is ordered as $\mathbb{N}$.
Moreover, each subword of $w$, which is ordered as $\mathbb{N}$ has this form.}

\medskip

If (2) were not true, then $W_k$ and $W_l$ would have a common infinite
subword. Then, by the Claim, we would have
$[p_{k(n)},q_{k(n)}][p_{k(n+1)},q_{k(n+1)}]\dots \equiv [p_{l(m)},q_{l(m)}][p_{l(m+1)},q_{l(m+1)}]\dots $
for some natural $n,m$. But this is impossible since $k,l$ are different in $\mathcal{K}$.


\medskip

c) Note that the group $H_1(Y)$ is torsion-free by \cite[Theorem~1.2]{Eda3} and \cite[Corollary~2.2]{Eda2}.
Therefore its minimal divisible subgroup containing $\{\overline{W_k}\,|\, k\in \mathcal{K}\}$ is isomorphic to the direct sum $\oplus_{i\in I}\mathbb{Q}$ for some index set $I$. Since the set $\{\overline{W_k}\,|\, k\in \mathcal{K}\}$ has cardinality $2^{\aleph_0}$, the set $I$ has cardinality $2^{\aleph_0}$ too.



\medskip

d) The statement of Theorem~9.2 concerning $H_1(Z)$ can be proven analogously, if
in~a) we take $a_i=[p_{2i-1},p_{2i}]$ and in b) we take $a_i=[p_{k(2i-1)},p_{k(2i)}]$.

\medskip

e) Now we prove, that $\pi_1(Y)$ contains a subgroup isomorphic to $\mathbb{Q}$.
In a) we have constructed a nontrivial element $W$ in $G=\mathcal{W}(P\cup Q)$ and noticed that there is a
natural homomorphism $\varphi: G\rightarrow \pi_1(Y)$.
By Theorem~8.1, $\varphi(W)$ has infinite order in $\pi_1(Y)$.

For any $i=2,3,\dots $, let $W_{i}$ denote the word obtained from $W$ by deletion of all occurrences
of $a_1,a_2,\dots ,a_{i-1}$. By Theorem 8.1, $\varphi(W_{i})=\varphi(W)$. We define a countable subword $U_i$ of $W_i$ as the subword, which starts from the first letter of $W_i$ (that is from $a_{i}$) and contains all letters of $W_{i}$ before
the second occurrence of $a_{i}$.

One can easily see, that $W_{i}\equiv (U_i)^{i!}$ and that $\varphi(U_{i})=\varphi((U_{i+1})^{i+1})$. Therefore the elements $\varphi(U_{i})$, $i=2,3,\dots$, generate a group with the presentation $\langle u_i\,\, (i=2,3,\dots)\,|\, u_i=u_{i+1}^{i+1}\rangle$. The group $\mathbb{Q}$ also has this presentation, as can be seen, if we interpret $u_i$ as $\frac{1}{i!}$.
\hfill $\Box$

\section{Acknowledgements}

The first named author thanks the MPIM at Bonn for its support and excellent working conditions
during the fall 2008, while this research was finished. He also thanks the University of Gdansk
for the hospitality during several visits in 2005-2008. The second named author acknowledges
financial support from the University of Gda\'nsk (BW UG 5100-5-0096-5) and from the Polish Ministry of Science (KBN 0524/H03/2006/31).

\end{document}